\pdfoutput=1 
\documentclass[11pt,a4paper]{article}

\usepackage[margin=30mm]{geometry}
\usepackage[T1]{fontenc}
\usepackage[utf8]{inputenc}
\usepackage{amsmath,amscd}
\usepackage{amsfonts}
\usepackage{amssymb}
\usepackage{amsthm}
\usepackage{xcolor}
\usepackage{enumitem}
\usepackage{tikz-cd}
\usepackage{pdfpages}
\usepackage{hyperref}
\DeclareMathOperator{\Pic}{Pic}

\newcommand{\PP}{\mathbb{P}}

\newcommand{\RR}{\mathbb{R}}

\newcommand{\ZZ}{\mathbb{Z}}

\newcommand{\CC}{\mathbb{C}}

\newtheorem{theorem}{Theorem}
\newtheorem{corollary}[theorem]{Corollary}

\newtheorem{lemma}[theorem]{Lemma}
\newtheorem{proposition}[theorem]{Proposition}

\newcommand{\C}{\CC}
\newcommand{\R}{\RR}
\newcommand{\Z}{\ZZ}
\newcommand{\F}{\mathbb{F}}
\newcommand{\OO}{\mathcal O}
\newcommand{\calC}{\mathcal C}
\newcommand{\ii}{\mathrm i}
\DeclareMathOperator{\wt}{wt}

\DeclareMathOperator{\Span}{span}
\DeclareMathOperator{\Sing}{Sing}
\DeclareMathOperator{\Sp}{Sp}

\title{Non-existence of quintic surfaces with 20 cusps and
of canonical maps of degree 9 in Beauville's second case}
\author{Lev Borisov \and Carlos Rito \and Xavier Roulleau}
\date{}

\begin{document}
\maketitle

\begin{abstract}
We prove that no complex quintic surface has exactly 20 ordinary cusps and no other singularities.
We also show
that a degree-9 canonical map onto a surface of positive geometric
genus would have such a quintic as its image. Hence the canonical degree
in Beauville's case II is at most 8. The proof combines the ternary
cusp code, abelian covers, Castelnuovo--de Franchis theorem and a
Hodge-theoretic bound on pairs of holomorphic 1-forms.
\end{abstract}

\section{Introduction}

Determining the maximum number of singularities on a surface of fixed
degree is a classical problem in algebraic geometry. Its difficulty is
already apparent for quintics, the hypersurfaces of smallest degree in
$\PP^3$ whose smooth members are of general type. Beauville
\cite{Beauville1980} determined the maximum number of ordinary nodes on
a quintic to be 31, using binary codes to study the relations
among exceptional curves. For ordinary cusps, that is, singularities
of type $A_2$, Varchenko's spectral bound \cite{Varchenko1983} gives an
upper bound of 20. Whether a quintic with 20 ordinary cusps
and no other singularities exists has remained open for more than
40 years.

The extremal cusp configuration is also tied to a fundamental question
about canonical maps. In his seminal 1979 paper \cite{Beauville1979},
Beauville proved that, if the canonical image $Y$ of a smooth surface
$T$ of general type is a surface, then either $p_g(Y)=0$, or $Y$ is
itself a canonical surface with $p_g(Y)=p_g(T)$. In the latter
situation, known as Beauville's case II for the degree of the canonical
map $\phi_{|K_T|}$, his inequalities give
\[
 \deg\varphi_{|K_T|}\le9.
\]
The possibility of equality has remained unresolved since that work.

In this paper, we show that equality would force the canonical image to be a quintic
with exactly 20 ordinary cusps, and then exclude the existence of such a quintic.
Consequently, if the canonical map of a smooth complex projective
surface of general type has a 2-dimensional image of positive
geometric genus, its degree is at most 8.

A geometric route to the extremal case was proposed by Tan
\cite[Theorem~I(4)]{Tan}. If the quotient of a fake projective plane by an
action of the group $(\Z/3)^2$ were a Godeaux surface with 4 ordinary cusps, its
\'etale cover of degree 5 would be a quintic with 20 cusps.
Base change would then produce a smooth surface whose canonical map
has degree 9 onto the quintic. The subsequent classification of
fake projective planes \cite{CartwrightSteger,PrasadYeung,PrasadYeungAddendum} revealed
examples with $(\Z/3)^2$-actions whose quotients are numerical Godeaux
surfaces with 4 cusps \cite[Corollary~1.2]{Keum2008}. Their
fundamental groups, however, are trivial or $\Z/2$, which prevents
the degree-5 \'etale covering required by Tan's construction
\cite{CartwrightSteger}. The classification therefore closed this
route to 20 cusps, while leaving the general existence question
open.

The central idea of our proof is to compare holomorphic forms arising
from different divisible sets of cusps. The cover associated with the
entire cusp code places these forms on a single surface. Their
products then impose a compatibility condition that is invisible in
the parameters and weight distribution of the code.

We first establish the link with canonical maps in
Section~\ref{sec:canonical-nine}. If a canonical map in Beauville's
case II has degree 9, then, after passing to the minimal model,
the numerical inequalities force
\[
 p_g(T)=4,\qquad q(T)=0,\qquad K_T^2=45
\]
(cf. \cite{Tan}).
The equality case of the Miyaoka--Yau theorem and the geometry of the
canonical system show that the map is finite, that its image is a
normal quintic with Du Val singularities, and that it is \'etale over
the smooth locus of the image. The degree restricts the local
quotient groups, leaving only singularities of types $A_2$ and $A_8$.
An Euler-characteristic calculation, combined with Varchenko's
semicontinuity theorem for the spectrum, excludes the $A_8$ points
and forces exactly 20 ordinary cusps.

To rule out a quintic with exactly twenty ordinary cusps and no
other singularities, we begin in Section~\ref{sec:quintics} with its ternary cusp code $\calC\subset\F_3^{20}$,
which records the relations
modulo 3 among the exceptional curves on the minimal resolution.
The discriminant of the intersection pairing over $\F_3$, together
with the Griesmer bound and Tan's restriction of the nonzero weights
to $12$, $15$, and $18$, forces
$\dim_{\F_3}\calC=5$. The MacWilliams identity and the minimum-weight
condition then give exactly 10 pairs $\{r_i,-r_i\}$ of
weight-18 words, with
\[
 \wt(r_i+r_j)=\wt(r_i-r_j)=12\qquad(i\ne j).
\]
These constraints are consistent: up to equivalence, there are exactly
2 ternary codes with the resulting parameters \cite{HHY}. The remaining issue is
whether such a code can be realized by the cusps of a quintic.

We take the normalized fibre product, over the minimal resolution of $Y$,
of the 5 cyclic triple covers associated with a basis of $\calC$.
Contracting the curves above the cusps gives a smooth surface $X$ and
a cover $X\to Y$ with group $G\simeq(\Z/3)^5$.
Every nonzero word $c\in \calC$ determines a character of $G$, and the
quotient of $X$ by the kernel of that character is the cyclic triple
cover of $Y$ associated with $c$.
For each of the 10
weight-18 pairs, Tan's theorem gives a 2-dimensional space of
holomorphic 1-forms on the resolved cyclic quotient. Choose a basis
of eigenforms. The Castelnuovo--de Franchis theorem and Beauville's
inequality for surface fibrations show that their wedge is nonzero.
All holomorphic 2-forms on this quotient are invariant, while no
nonzero holomorphic 1-form is invariant. Hence the 2 eigenforms have
opposite characters and pull back to forms $\alpha_i,\beta_i$ on $X$.
The mixed-weight identities imply that the $(1,1)$-eigenspaces for
$\chi_i+\chi_j$ and $\chi_i-\chi_j$ have dimension 1. Meanwhile,
the invariant real space $V=H^{1,1}(X,\R)^G$
has dimension 5 and signature $(1,4)$.

The decisive step is a general bound for finite abelian group actions
on smooth projective surfaces (Theorem~\ref{thm:pair}). If $m$ pairs
of holomorphic 1-forms have opposite characters and nonzero wedge
products, and the $(1,1)$-eigenspaces for sums and differences of
characters from distinct pairs have dimension at most 1, then
\[
 m\le\dim_{\R}H^{1,1}(X,\R)^G-1.
\]
Each pair determines a vector of negative square in the abstract real
exterior square of the invariant $(1,1)$ space, equipped with the form
induced by intersection. The mixed eigenspace conditions make these
vectors pairwise orthogonal. The Hodge index theorem controls the
negative index of the induced form and gives the bound. In the quintic
case it allows at most 4 such pairs, whereas the cusp code forces
10. This excludes their simultaneous realization on the full cover.

\medskip
\noindent\textbf{Use of AI tools.}
ChatGPT-6 (OpenAI) was used in this research project
to facilitate some technical arguments, namely in the proof of Theorem
7. The authors take full responsibility for the correctness of the
results of this paper.

\medskip
\noindent\textbf{Acknowledgments.}
The second author was financed by Portuguese funds through FCT
(Funda\c{c}\~ao para a Ci\^encia e a Tecnologia). This work was carried
out within project UID/00013/2025 of the Centro de Matem\'atica da
Universidade do Minho (CMAT/UM).

\medskip
\noindent\textbf{Notations and conventions.}
All varieties are complex and projective.
The symbol $\equiv$ denotes linear equivalence of divisors.
For an integral surface $Y$, the symbols $p_g(Y)$ and $q(Y)$ denote
the geometric genus and irregularity of a desingularization of $Y$.
We use additive notation for characters, so complex conjugation sends
$\chi$ to $-\chi$. A subscript $\chi$ denotes the corresponding
complex eigenspace, and
$H^{1,1}(X,\R)=H^{1,1}(X)\cap H^2(X,\R)$.
Signatures list the positive index first and the negative index second.

\section{Canonical maps of degree 9}\label{sec:canonical-nine}

We first show that a degree-9 canonical map in Beauville's case II
forces its image to be a quintic with exactly 20 ordinary cusps.

\begin{theorem}\label{thm:degree-nine}
Let $T$ be a smooth minimal surface of general type, and assume that its
canonical map
\[
  \varphi=\varphi_{|K_T|}\colon T \dashrightarrow
  \PP^{p_g(T)-1}
\]
is of degree $9$ onto a surface $Y$, with $p_g(Y)>0$.
Then $Y$ is a quintic surface in $\PP^3$ whose singular locus
consists of exactly 20 ordinary cusps.
\end{theorem}

\begin{proof}
\noindent\textit{Step 1: the numerical equality case.}
By Beauville's theorem \cite[Th\'eor\`eme~3.1]{Beauville1979}, since
$p_g(Y)>0$, the surface
$Y$ is of general type, is a canonical surface in Beauville's sense, and
\[
  p_g(Y)=p_g(T).
\]
Resolve the canonical map to a morphism $\widehat\varphi:\widehat T\to Y$,
and let $M=\widehat\varphi^*\OO_Y(1)$ be its free moving part.
Castelnuovo's inequality for a canonical surface
\cite[Remarque~5.6(1)]{Beauville1979} gives
$\deg Y\geq 3p_g(T)-7$. Hence
\begin{equation}\label{eq:canonical-chain}
  9(3p_g(T)-7)
  \leq 9\deg Y
  =M^2
  \leq K_T^2
  \leq 9\chi(\OO_T)
  =9\bigl(p_g(T)-q(T)+1\bigr).
\end{equation}
The last inequality is the Bogomolov--Miyaoka--Yau inequality.
It follows that
\[
  2p_g(T)+q(T)\leq 8.
\]
Since the canonical image is a surface in $\PP^{p_g(T)-1}$ and is of general
type, one has $p_g(T)\geq4$: for $p_g(T)=3$ its image would be all of $\PP^2$.
Therefore
\[
  p_g(T)=4,\qquad q(T)=0.
\]
Every inequality in \eqref{eq:canonical-chain} is now an equality.  Thus
\[
  \deg Y=5,\qquad K_T^2=M^2=45.
\]
Since $\chi(\OO_T)=1+p_g(T)-q(T)=5$, Noether's formula yields
\[
  c_2(T)=12\chi(\OO_T)-K_T^2=60-45=15.
\]
Thus $K_T^2=3c_2(T)$. By the equality case of the Miyaoka--Yau theorem,
the canonical model of $T$ is smooth \cite[Theorem~1.2]{GKPT}.
Since $T$ is a minimal surface, it coincides with this model.
In particular, $K_T$ is ample.

\medskip
\noindent\textit{Step 2: the canonical system has neither a fixed part nor
base points.}
Write
\[
  |K_T|=F+|L|,
\]
where $F$ is the fixed divisor and $|L|$ has no fixed components.  If the
base points of $|L|$ (possibly infinitely near) have multiplicities $m_i$, then
\[
  M^2=L^2-\sum_i m_i^2.
\]
Since $M^2=K_T^2$, one has $L^2\geq K_T^2$.  On the other hand, if
$F\neq0$, then
\[
  K_T^2-L^2=(K_T-L)(K_T+L)=F\cdot(K_T+L)>0.
\]
Indeed, $K_T\cdot F>0$ because $K_T$ is ample, while $L\cdot F\geq0$
because a general member of $|L|$ has no common component with $F$.
This is a contradiction.  Hence $F=0$, and then
\[
  M^2=K_T^2-\sum_i m_i^2=K_T^2
\]
forces every $m_i$ to vanish.  Therefore $|K_T|$ is base-point-free and
$  K_T=\varphi^*\OO_Y(1).$ 
Since $K_T$ is ample, $\varphi$ cannot contract a curve.  Thus it is a
finite morphism of degree $9$.

\medskip
\noindent\textit{Step 3: normality, canonical singularities, and
quasi-\'etaleness.}
Let $\nu:Y^\nu\to Y$ be the normalization of $Y$, and let
$\mu:S\to Y^\nu$ be its minimal resolution. Set
\[
  \eta=\nu\circ\mu:S\longrightarrow Y,
  \qquad H=\eta^*\OO_Y(1).
\]
By Beauville's theorem \cite[Th\'eor\`eme~3.1, cas II]{Beauville1979}, since we are in the case where the canonical image
is of general type with positive geometric genus, the canonical map of
$S$ is the composition
\[
  S\xrightarrow{\eta}Y\hookrightarrow\PP^3.
\]
Consequently $|H|$ is the moving part of the canonical system of $S$, so
\begin{equation}\label{eq:canonical-system}
  |K_S|=E+|H|,
  \qquad E\geq0,
\end{equation}
where $E$ is the fixed part. See \cite[Th\'eor\`eme~3.1 and \S3.3]{Beauville1979}.

Choose a general member $C_S\in|H|$.  By Bertini's theorem it is smooth
and irreducible, and it maps birationally onto a general plane section
$  C=Y\cap\PP^2$, 
which is an integral plane quintic.  Moreover, one has $  H^2=\deg Y=5.$ 
Since $C_S$ is smooth, adjunction on $S$, together with
\eqref{eq:canonical-system}, gives
\[
\begin{aligned}
  2g(C_S)-2
    &=(K_S+H)\cdot H=(2H+E)\cdot H=10+E\cdot H.
\end{aligned}
\]
The divisor $H$ is nef and $E$ is effective, so $E\cdot H\geq0$ and hence
$g(C_S)\geq6$.  On the other hand, $C_S$ is the normalization of the
integral plane quintic $C$, whose arithmetic genus is $6$.  Therefore $ g(C_S)\leq6.$
 It follows that
\begin{equation}\label{eq:canonical-EH}
  E\cdot H=0,
  \qquad
  g(C_S)=6,
\end{equation}
and that the general plane section $C$ is smooth.
Therefore $Y$ is normal.
Thus $\nu$ is an isomorphism, and $\eta:S\to Y$ is the minimal
resolution of $Y$.

Adjunction now gives
\[
  K_Y=\OO_Y(1),
  \qquad
  H=\eta^*K_Y.
\]
By \eqref{eq:canonical-EH}, every irreducible component of $E$ is contracted by
$\eta$, because $H$ is the pullback of an ample divisor. We have then
\[
  K_S=\eta^*K_Y+E.
\]
The discrepancy divisor $E$ is effective, hence $Y$ has canonical
singularities only, and then $E=0$. In dimension 2 these are precisely
the Du Val singularities. In particular, they are rational. This implies $ q(Y)=q(S)=0.$

Finally, the ramification formula for the finite morphism $\varphi$ reads
\[
  K_T=\varphi^*K_Y+R_\varphi=\varphi^*\OO_Y(1)+R_\varphi.
\]
Comparison with $K_T=\varphi^*\OO_Y(1)$ shows that $R_\varphi=0$.
Thus $\varphi$ is \'etale over $Y_{\mathrm{reg}}$, and hence $\varphi$
is quasi-\'etale.

\medskip
\noindent\textit{Step 4: the possible local singularities.}
Fix $y\in\Sing(Y)$.  Since $Y$ has a Du Val singularity at $y$, one has
\[
  (Y,y)\cong(\C^2/G_y,0)
\]
for a finite subgroup $G_y\subset\operatorname{SL}_2(\C)$.
Since $T$ is smooth and $\varphi$ is quasi-\'etale, the map near each
$x\in\varphi^{-1}(y)$ is analytically equivalent to the quotient map
$\C^2\to\C^2/G_y$. Thus every local degree over $y$ equals $|G_y|$.
These local degrees sum to $9$, so $|G_y|$ divides $9$.
The only nontrivial finite subgroups of $\operatorname{SL}_2(\C)$ whose
orders divide $9$ are the cyclic groups of orders $3$ and $9$.  They give,
respectively,
\[
  A_2=\frac13(1,2),
  \qquad
  A_8=\frac19(1,8).
\]
Thus the singularities of $Y$ are only of types $A_2$ and $A_8$.  Write
\[
  a=\#\{A_2\text{-points of }Y\},
  \qquad
  b=\#\{A_8\text{-points of }Y\}.
\]

\medskip
\noindent\textit{Step 5: the Euler-number equation.}
Since $\eta:S\to Y$ is the minimal resolution of the Du Val
singularities of $Y$, then
\[
  K_S^2=5,
  \qquad
  p_g(S)=4,
  \qquad
  q(S)=0.
\]
Noether's formula gives
\[
  e(S)=12\chi(\OO_S)-K_S^2=60-5=55.
\]
The exceptional divisor over an $A_2$-point has Euler number $3$, and the
exceptional divisor over an $A_8$-point has Euler number $9$.
Let $D\subset S$ be the exceptional locus of $\eta$. The resolution
restricts to an isomorphism $S\setminus D\cong Y_{\mathrm{reg}}$.
By additivity of the Euler characteristic,
\[
  e(Y_{\mathrm{reg}})=e(S)-e(D)=55-3a-9b.
\]
There are $9/3=3$ points of $T$ over an $A_2$-point and 1 point over an
$A_8$-point.  Since $\varphi$ is an unramified covering of degree $9$ over
$Y_{\mathrm{reg}}$, one obtains
\[
\begin{aligned}
  15=e(T)
  &=9e(Y_{\mathrm{reg}})+3a+b=9(55-3a-9b)+3a+b.
\end{aligned}
\]
Equivalently,
\begin{equation}\label{eq:canonical-basket}
  3a+10b=60.
\end{equation}

\medskip
\noindent\textit{Step 6: excluding $A_8$-points.}
Varchenko's theorem bounds the total number of local spectral numbers
in an open interval of length 1 by the corresponding number for
$x^5+y^5+z^5$. See \cite[\S\S1--2]{Varchenko1983} and
\cite[\S5.4]{vanStraten2020}.
(The statement in \cite[\S1.5]{vanStraten2020} contains an indexing typo:
for a hypersurface in $\PP^n$, the comparison polynomial should read
$x_1^d+\cdots+x_n^d$, with $n$ variables.)
In the unshifted convention, the spectrum of $x^5+y^5+z^5$ is the multiset
\[
\Sp(x^5+y^5+z^5)=
\left\{\frac{i+j+k}{5}:1\le i,j,k\le4\right\}.
\]
The interval $I=(6/5,11/5)$ contains the terms with
$i+j+k=7,8,9,10$, numbering $12+12+10+6=40$.
The spectrum of an $A_n$-singularity is
$\{1+j/(n+1):1\le j\le n\}$, so an $A_2$-point contributes 2
numbers to $I$ and an $A_8$-point contributes 7. Hence
\[
2a+7b\le40.
\]
Together with $3a+10b=60$ from \eqref{eq:canonical-basket}, this gives
$120+b\le120$. Thus $b=0$ and $a=20$.
\end{proof}

\section{Non-existence of 20-cuspidal quintics}\label{sec:quintics}

\subsection{The exceptional lattice forces a 5-dimensional code}\label{sec:lattice}

Let $Y\subset\PP^3$ be a quintic surface with exactly 20
cusps and no other singularities, and let $\eta:S\to Y$ be its smooth
minimal resolution. 
Write $H\in H^2(S,\Z)$ for the class of $\eta^*\OO_Y(1)$.
One has $K_S=H$ and the invariants of $S$ are 
\begin{equation}\label{eq:invariants}
 K_S^2=5,\quad p_g(S)=4,\quad q(S)=0,
\quad\chi(\OO_S)=5,\quad e(S)=55.
\end{equation}
Let $E_i,E'_i$ be the exceptional curves over the $i$-th cusp.
Following Barth--Rams \cite[Introduction]{BRcodes}, the
\emph{ternary cusp code of $Y$} is $\calC=\calC(Y):=\ker\varphi$, where
\begin{equation}\label{eq:cuspcode}
\varphi:\F_3^{20}\longrightarrow H^2(S,\F_3),\qquad
(c_1,\ldots,c_{20})\longmapsto
\sum_{i=1}^{20}c_i\bigl([E_i]-[E'_i]\bigr).
\end{equation}
Here the brackets denote cohomology classes modulo 3.
Elements of $\mathcal C$ are called words.
For $c=(c_1,\ldots,c_{20})\in\mathcal C$, its weight is
$  \operatorname{wt}(c)=\#\{i:c_i\ne0\}$. 
The Hamming distance between 2 codewords $c,d$ is
$\operatorname{wt}(c-d)$. The length of the code $\mathcal C$, i.e. the dimension of its ambient space, is $20$. 
Define $k=\dim_{\F_3}\calC$.
Let us prove $k=5$ from the intersection pairing modulo 3.

\begin{lemma}\label{lem:ambientcohomology}
The space $P=H^2(S,\F_3)$ has dimension $53$, and its intersection
pairing is nondegenerate with Gram determinant $1$.
\end{lemma}
\begin{proof}
Apply the Brieskorn--Tjurina simultaneous resolution theorem
\cite[Theorem~168]{Catanese} after a finite base change to a smoothing of $Y$,
and then Ehresmann's theorem \cite[\S12.2]{Catanese}.
This shows that $S$ is diffeomorphic to a smooth quintic.
Thus $S$ is simply connected. Consequently, $\Lambda=H^2(S,\Z)$ is
torsion-free.
By integral Poincar\'e duality, every integral basis has an integral
dual basis for the intersection pairing. Thus its Gram matrix $B$
has an integral inverse, so $\det B=\pm1$.

Using $b_1(S)=b_3(S)=0$, $e(S)=55$, and $p_g(S)=4$, we obtain
$b_2(S)=53$ and $h^{1,1}(S)=45$.
By the Hodge index theorem \cite[Corollary~3.3.16]{Huybrechts}, the
intersection form on $H^2(S,\R)$ has signature
$\bigl(1+2p_g(S),\,h^{1,1}(S)-1\bigr)=(9,44)$.
Since the negative index is even, $\det B$ is positive and hence equals $1$.

Finally, Poincar\'e duality gives $H^3(S,\Z)\cong H_1(S,\Z)=0$.
The coefficient sequence for reduction modulo 3 therefore identifies
$P$ with $\Lambda/3\Lambda$. Reducing an integral basis and its Gram
matrix gives $\dim P=53$ and $\det P=1$ in $\F_3$.
Over $\F_3$, a change of basis preserves a nonzero Gram determinant,
since every nonzero element has square $1$.
\end{proof}

\begin{proposition}\label{prop:kfive}
The cusp code has dimension $k=5$.
\end{proposition}
\begin{proof}
For a quintic with 20 cusps, Borisov--Rito
\cite[Section~2, proof of Proposition~1]{BoRi} give $k\ge4$.
Tan's result \cite[Theorem~4.3.1]{Tan}, recalled in
\cite[Introduction]{BRcodes}, restricts the nonzero weights to
$\{12,15,18\}$, so the minimum distance between 2 words in $\mathcal C$ is at least 12.
Suppose that $k\ge6$. A 6-dimensional subcode $\mathcal C _6$ would have the same lower bound
on its minimum distance. 
For a ternary linear code of dimension $k$, minimum distance $d$,
and length $n$, the Griesmer bound \cite[Corollary~2.1.4]{Tan}  states that
$  n\geq\sum_{j=0}^{k-1}\left\lceil\frac{d}{3^j}\right\rceil.$ 
Applying it to $\mathcal C _6$ gives
\[
  n\geq\sum_{j=0}^{5}\left\lceil\frac{12}{3^j}\right\rceil
    =12+4+2+1+1+1=21,
\]
whereas $\mathcal C _6$ has length $n=20$. Hence $4\le k\le5$.

Suppose, for a contradiction, that $k=4$.

\smallskip
\emph{Step 1: an orthogonal complement of dimension $32$ and determinant $2$.}
Each exceptional pair $E_i,E'_i$ satisfies $E_i^2={E'_i}^2=-2$,
$E_i E'_i=1$ and $H\cdot E_i=H\cdot E'_i=0$.
Write the same letters for the reductions of these classes in $P$, and put
\[
s_i=E_i+E'_i,\qquad v_i=E_i-E'_i\qquad(1\le i\le20).
\]
The intersection matrix gives, modulo 3,
\[
s_i^2=1,\qquad v_i\cdot v_j=0,\qquad s_i\cdot v_j=0
\quad\text{for all }i,j,
\]
while $H^2=2\pmod3$. 

The vectors $H,s_1,\ldots,s_{20}$ have Gram matrix
$\operatorname{diag}(2,1,\ldots,1)$ over $\F_3$, which is invertible. 
The vectors span therefore a $21$-dimensional nondegenerate subspace
$ A:=\Span_{\F_3}(H,s_1,\ldots,s_{20})\subset P.$
Since $A$ is nondegenerate, $P=A\oplus A^\perp$.
By Lemma~\ref{lem:ambientcohomology}, $\dim P=53$ and $\det P=1$. Hence
\begin{equation}\label{eq:mod3complement}
\dim A^\perp=53-21=32,\qquad \det(A^\perp)=\frac{\det P}{\det A}=\frac12=2
\quad\text{in }\F_3.
\end{equation}

\smallskip
\emph{Step 2: the code gives a 16-dimensional totally isotropic subspace.}
For the map \eqref{eq:cuspcode}, put
$U:=\operatorname{im}\varphi=\Span_{\F_3}(v_1,\ldots,v_{20})$.
We have
\begin{equation}\label{eq:codekernel}
\dim U=20-\dim\ker\varphi=20-k=16.
\end{equation}
Each $v_i$ is orthogonal to $H$ and all $s_j$, so $U\subset A^\perp$.
Moreover, $v_i\cdot v_j=0$ for every $i,j$, so the pairing is totally isotropic on $U$: it vanishes
identically on $U\times U$.

\smallskip
\emph{Step 3: this subspace forces determinant $1$.}
Extend a basis of $U$ to a basis of $A^\perp$. Since $U$ is totally isotropic,
the Gram matrix in this basis has the form
\[
\begin{pmatrix}
0 & C \\
C^{\mathsf T} & B
\end{pmatrix}
\]
with $16\times16$ blocks. Since $16$ is even, interchanging the 2
blocks of columns leaves the determinant unchanged. Thus
$\det(A^\perp)=(\det C)^2$, contradicting \eqref{eq:mod3complement} because
$2$ is not a square in $\F_3$. Hence $k=5$.
\end{proof}

\subsection{The 10 opposite pairs in the cusp code}\label{sec:code}

We thus know that the code $\mathcal C$ is $5$-dimensional.

\begin{proposition}\label{prop:enumerator}
The code $\calC$ has full support and exactly 10 pairs of opposite
weight-18 words. If $r_1,\ldots,r_{10}$ represent these pairs, then
\begin{equation}\label{eq:weightsmixed}
\wt(r_i+r_j)=\wt(r_i-r_j)=12\qquad(i\ne j).
\end{equation}
\end{proposition}
\begin{proof}
If a coordinate vanished identically, puncturing it would give a ternary
code of length 19, dimension 5 and minimum distance at least
12. The Griesmer bound  instead requires length at least
$12+4+2+1+1=20$ (see Proof of Proposition \ref{prop:kfive}). 
Thus every coordinate is nonzero on some word: the code has full support.

Let $A_w$ count codewords of weight $w$, and set $a=A_{18}$.
At each coordinate, exactly $2\cdot3^4=162$ words are nonzero.
Count the pairs $(c,i)\in\calC\times\{1,\ldots,20\}$ such that $c_i\ne0$.
Counting first by codeword gives the sum of their weights; counting first
by coordinate gives $20\cdot162$. Thus
\[
12A_{12}+15A_{15}+18a=20\cdot162=3240.
\]
There are
$3^5-1=242$ nonzero words. Hence $A_{12}+A_{15}+a=242$.
These 2 equations give $A_{12}=130+a$ and $A_{15}=112-2a$.

For a code $D$ of length $20$, write
$W_D(x,y)=\sum_{d\in D}x^{20-\wt(d)}y^{\wt(d)}$ for its weight enumerator.
Thus
\[
W_{\calC}(x,y)=x^{20}+(130+a)x^8y^{12}
               +(112-2a)x^5y^{15}+ax^2y^{18}.
\]
Let $B_w$ count the words of weight $w$ in the dual code
$$\calC^\perp=\{u\in\F_3^{20}:u\cdot c=0\text{ for every }c\in\calC\}.$$
Here $u\cdot c = \sum_1^{20}u_ic_i$
denotes the standard dot product on $\mathbb F_3^{20}$.
The dual code $\calC^\perp$ is the orthogonal complement of $\calC$ with respect
to this form.
By the MacWilliams identity~\cite[Chap.~5, \S~6, Thm.~13, Eq.~(47), p.~146]{MS},
\[
W_{\calC^\perp}(x,y)=3^{-5}W_{\calC}(x+2y,x-y).
\]
Comparing the coefficients of $x^{18}y^2$ and $x^{17}y^3$ gives
\begin{equation}\label{eq:MacWilliams}
B_2=\frac{a-20}{3},\qquad B_3=2(20-a).
\end{equation}
Since both are nonnegative, $a=20$. The words of weight 18
therefore form 10 pairs $\{r_i,-r_i\}$.

For $i\ne j$, let $s$ be the number of coordinates where both $r_i$
and $r_j$ are nonzero. Since both words have weight 18,
$s\ge 18+18-20=16$. At each of these $s$ coordinates, exactly 1 of
$r_i+r_j$ and $r_i-r_j$ is nonzero over $\F_3$. Thus
\[
\wt(r_i+r_j)+\wt(r_i-r_j)\le s+2(20-s)=40-s\le24.
\]
The 2 representatives are not proportional, so their sum and difference
are nonzero codewords. Each has weight at least 12, and hence both
have weight 12, as claimed in \eqref{eq:weightsmixed}.
\end{proof}
These numerical constraints do not yet give a contradiction. Indeed,
up to monomial equivalence, there are exactly 2 ternary
$[20,5,12]$ codes, both with the weight distribution found above
(see \cite{HHY}). Below we use a geometric argument to
show that such a code cannot arise from 20 cusps on a quintic.

\subsection{The \texorpdfstring{$(\Z/3)^5$}{(Z/3)\string^5}-cover and its character spaces}\label{sec:covers}

\subsubsection*{Construction and smoothness}
We construct a smooth $(\Z/3)^5$-cover $X\to Y$ whose cyclic triple
quotients are the covers associated with the nonzero words of $\cal C$.
We first take the normalized
fibre product of 5 cyclic triple covers over the minimal resolution
$\eta:S\to Y$. Write $p_i$ for the cusp corresponding to $E_i,E'_i$ and
$[x]$ for the representative of $x\in\F_3$ in $\{0,1,2\}$.
For a non-zero word $c=(c_i)\in\calC$, the code relation \eqref{eq:cuspcode}
gives a divisor $L_c$ on $S$ with
\[
  3L_c\equiv B_c:=\sum_{i=1}^{20}
    \bigl([c_i]E_i+[-c_i]E'_i\bigr).
\]
Its class is unique since $\Pic(S)$ has no $3$-torsion. Let $V_c\to S$ be the normal cyclic triple cover associated with
$3L_c\equiv B_c$. The Stein factorization of $V_c\to S\to Y$ gives
a commutative diagram
\[
\begin{tikzcd}[column sep=3.6em,row sep=1.3em]
 V_c \arrow[r] \arrow[d] & T_c \arrow[d,"\pi_c"]\\
 S \arrow[r,"\eta"'] & Y,
\end{tikzcd}
\]
where the finite map $\pi_c$ ramifies precisely at the $p_i$ with $c_i\ne0$,
and $T_c$ is smooth above those cusps \cite[Lemma~3.1.1]{Tan}.

By Propositions~\ref{prop:kfive} and~\ref{prop:enumerator}, $\calC$ has
dimension 5 and full support. Choose any basis $c^{(1)},\ldots,c^{(5)}$
of $\calC$ and take the normalized fibre product of the corresponding
finite flat covers of $S$:
\[
  \widehat\pi:\widehat X=
   \operatorname{Norm}\bigl(V_{c^{(1)}}\times_S\cdots\times_S
    V_{c^{(5)}}\bigr)\longrightarrow S.
\]
Away from the exceptional curves, $\widehat{\pi}$ is \'etale, and
$G=(\Z/3)^5$ acts freely and transitively on each fibre. Let $M$ be the matrix whose rows are the chosen basis
words, and let $m_i=(c_i^{(1)},\ldots,c_i^{(5)})^{\mathsf T}$ be its
$i$-th column. The branch coefficients of the 5 cyclic factors along
$E_i$ show that the inertia subgroup there is generated by $m_i$. Since
$M$ has rank 5, the $m_i$ span $G$.

The group $G$ acts transitively on the connected components of $\widehat X$,
since each orbit of components gives an open and closed subset of the
connected quotient $S=\widehat X/G$. The stabilizers of the components
are conjugate, hence equal because $G$ is abelian. Each inertia subgroup
fixes a point above its branch curve and therefore preserves the component
containing that point. The common stabilizer thus contains all the $m_i$
and equals $G$. By transitivity, $\widehat X$ is connected
(see also \cite[\S2.1, p.~309]{AlexeevPardini}).

The Stein factorization of $\eta\circ\widehat\pi$ gives the diagram
\[
\begin{tikzcd}[column sep=3.6em,row sep=1.3em]
 \widehat X \arrow[r,"\rho"] \arrow[d,"\widehat\pi"'] & X \arrow[d,"\pi"]\\
 S \arrow[r,"\eta"'] & Y .
\end{tikzcd}
\]
Here $\rho$ contracts the connected configurations above the cusps and
$\pi$ is a finite $G$-cover, \'etale over $Y_{\mathrm{reg}}$.
Fix a primitive cube root of unity $\zeta$ and write
$\chi_a(g)=\zeta^{a\cdot g}$ for $a\in\F_3^5$. The character $\chi_a$
corresponds to the codeword $c(a)=\sum_{j=1}^{5}a_jc^{(j)}$, since
$\chi_a(m_i)=\zeta^{a\cdot m_i}=\zeta^{\,c(a)_i}$ for each $i$. Thus
$\chi_a\leftrightarrow c(a)$ identifies $G^*$ with the row space $\calC$
of $M$. 
We write the group of complex characters $G^*$ additively: for
$\chi,\psi\in G^*$, the character $\chi+\psi$ is their pointwise
product, and $-\chi$ is the inverse character.
The weight of a character is the number of nonzero coordinates
of its corresponding codeword. For $a\ne0$, the cyclic quotient has branch data
$B_{c(a)},L_{c(a)}$ on $S$, hence $X/\ker\chi_a\simeq T_{c(a)}$ over $Y$.

At each cusp $p_i$, full support gives a basis word $c$ with $c_i\ne0$.
The corresponding $T_c$ is smooth above $p_i$. The finite projection
$X\to T_c$ is \'etale away from the points over the cusps, so by
purity of the branch locus \cite[\href{https://stacks.math.columbia.edu/tag/0BMB}{Tag~0BMB}]{stacks-project}
it is \'etale above $p_i$ as well. Hence $X$ is smooth.

\subsubsection*{Holomorphic forms on cyclic quotients}
For a subgroup $K\subset G$, put $T=X/K$ and let $Z\to T$ be its minimal
resolution. Pullback along a resolution of the rational map
$X\dashrightarrow Z$ gives $G/K$-equivariant identifications
\begin{equation}\label{eq:formscomparison}
H^{p,0}(Z)\simeq H^{p,0}(X)^K\qquad(p=1,2).
\end{equation}
Since quotient singularities are rational and taking invariants is exact,
we have $H^j(Z,\OO_Z)\simeq H^j(X,\OO_X)^K$.
Functoriality of Hodge decomposition and complex conjugation give
\eqref{eq:formscomparison}. Pullback commutes with exterior products,
so these identifications preserve wedge products and whether a wedge
vanishes.
For $\chi\in G^*$, write
\[
H^{1,0}(X)_\chi
=\{\alpha\in H^{1,0}(X):
g^*\alpha=\chi(g)\alpha
\text{ for every }g\in G\}.
\]
For each $g\in G$, the operator $g^*$ is diagonalizable over $\C$:
it satisfies $(g^*)^3=1$, and $t^3-1$ has 3 distinct roots.
Since $G$ is abelian, the other generators preserve the eigenspaces
of the first. Splitting these spaces successively into eigenspaces
for the remaining generators gives a common eigenbasis. The eigenvalues
on each common eigenline define a character of $G$. Hence
\[
H^{1,0}(X)=\bigoplus_{\chi\in G^*}H^{1,0}(X)_\chi.
\]

\subsubsection*{The invariant and weight-12 \texorpdfstring{$(1,1)$}{(1,1)} spaces}

\begin{proposition}\label{prop:charspaces}
The invariant space $V=H^{1,1}(X,\R)^G$ has dimension 5 and signature $(1,4)$.
Every character $\chi\in G^*$ of weight 12 occurs with multiplicity 1
in the $G$-representation $H^{1,1}(X)$. Equivalently,
\begin{equation}\label{eq:char11}
\dim H^{1,1}(X)_\chi=1.
\end{equation}
\end{proposition}
\begin{proof}
The quotient and resolution maps identify $V$ with the orthogonal
complement of the 40 exceptional classes in $H^{1,1}(S,\R)$. Under
this identification, the pairing on $V$ is $|G|$ times the pairing
on $S$. The exceptional classes span a negative definite subspace of
dimension 40. Since $H^{1,1}(S,\R)$ has signature $(1,44)$ by the
Hodge index theorem, $V$ has signature $(1,4)$ and hence dimension 5.

Now fix $\chi$ of weight 12, put $K=\ker\chi$, and let $Z\to T=X/K$
be the minimal resolution.
Tan's theorem \cite[Theorem~I(3)(a)]{Tan}
gives $p_g(Z)=6$ and $q(Z)=0$.
 Moreover $K_Z^2=3K_Y^2=15$ by \cite[Section~3.1, formula~(4)]{Tan}.
Noether's formula and Hodge decomposition yield
\[
e(Z)=12(1+6)-15=69,\qquad b_2(Z)=67,\qquad
h^{1,1}(Z)=67-2\cdot6=55.
\]
The other 8 cusps of $Y$ are unramified in $T\to Y$, so each has
3 $A_2$ points above it. Thus $T$ has 24 cusps and $Z\to T$
has 48 independent exceptional $(1,1)$ classes. The same
quotient and resolution comparisons therefore give
\[
\dim_\C H^{1,1}(X)^K=55-48=7.
\]
The only characters of $G$ trivial on $K$ are $0,\chi,-\chi$, so
\[
H^{1,1}(X)^K=H^{1,1}(X)^G
\oplus H^{1,1}(X)_\chi\oplus H^{1,1}(X)_{-\chi}.
\]
The first summand is the complexification of $V$ and has dimension 5.
Complex conjugation preserves type $(1,1)$ and exchanges the other 2
summands, so their dimensions agree. Each therefore has dimension $1$.
\end{proof}

\subsection{Opposite 1-forms have nonzero wedge}\label{sec:wedges}

\begin{proposition}\label{prop:nonzero}
Let $\chi\in G^*$ have weight 18. Then
\[
 \dim_\C H^{1,0}(X)_\chi=\dim_\C H^{1,0}(X)_{-\chi}=1,
\]
and any nonzero forms
\[
 \alpha\in H^{1,0}(X)_\chi \quad\text{and}\quad
 \beta\in H^{1,0}(X)_{-\chi} \quad\text{satisfy}\quad
 \alpha\wedge\beta\ne0.
\]
\end{proposition}
\begin{proof}
Put $K=\ker\chi$ and let $Z\to T=X/K$ be the minimal resolution.
By Tan's theorem \cite[Theorem~I(3)(c)]{Tan}, $q(Z)=2$. Let $u,v$ be
any basis of $H^{1,0}(Z)$ and suppose that $u\wedge v=0$. By the
Castelnuovo--de Franchis theorem \cite{BHPV}, there is a morphism
$f:Z\to B$ with connected fibres onto a smooth curve of genus $b\ge2$.

Write $\psi:Z\to Y$ for the composite $Z\to T\to Y$ and $H_Y$ for the
hyperplane class of $Y$. The cyclic triple cover $T\to Y$ is \'etale in
codimension 1, and $T$ has only $A_2$ singularities, so its minimal
resolution is crepant. Since $K_Y\equiv H_Y$, we have
\begin{equation}\label{eq:Kcyclic}
K_Z\equiv\psi^*H_Y,\qquad K_Z^2=3K_Y^2=15.
\end{equation}
In particular, $K_Z$ is nef and $Z$ is minimal.

Each exceptional curve of $Z\to T$ is rational and hence maps to a
point of $B$. A general fibre $F$ of $f$ therefore avoids the exceptional
curves. Write $g$ for its genus, $C=\psi(F)$ for its image and
$d$ for the degree of $F\to C$. Since $F^2=0$, adjunction on $Z$ and
the projection formula give
\begin{equation}\label{eq:fibredegree}
2g-2=K_Z\cdot F=\psi^*H_Y\cdot F=d\,\deg(C)>0.
\end{equation}
If $g=2$, then \eqref{eq:fibredegree} gives $(d,\deg(C))=(1,2)$ or $(2,1)$.
The first case would make $F$ birational to an integral plane conic,
contrary to $g=2$. In the second case, $C$ would be a line and the map $F\to C$
would be generically 2-to-1. The image of $F$ in $T$ would then contain
exactly 2 of the 3 points above a general point of $C$. But the deck
group permutes the components of the inverse image of $C$ in $T$ and
acts transitively on these 3 points. Each component must therefore
contain the same number of them, either 1 or 3. Hence $g\ge3$.

Beauville's inequality \cite[Appendix, Corollary~(i), p.~344]{Beauville1982}
applied to $f:Z\to B$ now gives
\[
15=K_Z^2\ge8(g-1)(b-1)\ge16,
\]
a contradiction. Thus $u\wedge v\ne0$ for every basis $u,v$ of $H^{1,0}(Z)$.

Now let $\sigma$ generate the cyclic group $\Gamma=G/K$ of order 3.
Since $(\sigma^*)^3=1$, we may choose $u,v$ to be eigenvectors of
$\sigma^*$, with eigenvalues $\lambda,\mu$. Applying
\eqref{eq:formscomparison} to $K$ and to $G$ gives
$H^{p,0}(Z)^\Gamma\simeq H^{p,0}(S)$ for $p=1,2$.
Since $q(S)=0$, neither eigenvalue is $1$. Moreover,
$p_g(S)=p_g(Z)=4$ by \eqref{eq:invariants} and Tan's theorem, so every
holomorphic 2-form on $Z$ is $\Gamma$-invariant. In particular,
\[
 u\wedge v=\sigma^*(u\wedge v)=\lambda\mu\,u\wedge v.
\]
Since $u\wedge v\ne0$, we have $\lambda\mu=1$. The eigenvalues are
therefore $\zeta$ and $\zeta^2$ in some order. The corresponding
characters of $\Gamma$ pull back to $\chi$ and $-\chi$ on $G$, so
\eqref{eq:formscomparison} identifies their eigenlines with
$H^{1,0}(X)_\chi$ and $H^{1,0}(X)_{-\chi}$. The wedge of nonzero forms
from these lines is a nonzero multiple of the pullback of $u\wedge v$.
\end{proof}

\subsection{A character-pair bound}\label{sec:pair}
We now combine the nonzero wedge products of opposite-character 1-forms
with the fact that the relevant $(1,1)$-eigenspaces have dimension 1
to bound the number of such pairs. 
\begin{theorem}\label{thm:pair}
Let a finite abelian group $G$ act on a smooth connected projective complex
surface $X$, and put $V=H^{1,1}(X,\R)^G$.
Let  $\chi_1,\ldots,\chi_m$ be $m$ characters of  $G$ and choose $m$ pairs of holomorphic 1-forms
\[
\alpha_i\in H^{1,0}(X)_{\chi_i},\qquad
\beta_i\in H^{1,0}(X)_{-\chi_i},\qquad 1\le i\le m.
\]
Suppose that $\alpha_i\wedge\beta_i\ne0$ for every $i$, and that, for
$i\ne j$,
\begin{equation}\label{eq:mixedbound}
\dim_\C H^{1,1}(X)_{\chi_i-\chi_j}\le1,
\qquad
\dim_\C H^{1,1}(X)_{\chi_i+\chi_j}\le1.
\end{equation}
Then
\[
m\le\dim_\R V-1.
\]
\end{theorem}

\begin{proof}
Let $H$ be an ample class on $X$. Its average
$H_G=|G|^{-1}\sum_{g\in G}g^*H$ is ample because each $g^*H$ is ample,
and is $G$-invariant because $G$ permutes the summands. Thus $H_G\in V$.
The intersection form $Q(A,B)=\int_X A\cup B$ is positive on $\R H_G$
and, by the Hodge index theorem \cite{BHPV}, negative definite on its
orthogonal complement in $H^{1,1}(X,\R)$. Its restriction to $V$ therefore
has signature $(1,t)$, where $t=\dim_\R V-1$.

Set $U_i=\alpha_i\wedge\beta_i$ and
\[
A_i=\ii[\alpha_i\wedge\overline{\alpha_i}],\qquad
B_i=\ii[\beta_i\wedge\overline{\beta_i}].
\]
Since $X$ is compact K\"ahler, its holomorphic 1-forms are closed.
The classes $A_i$ and $B_i$ are real: conjugation changes the sign of
$\ii$ and reverses the order of the 1-forms in each wedge.
For $g\in G$, put $\lambda=\chi_i(g)$. Since $g^*\alpha_i=\lambda\alpha_i$
and pullback commutes with complex conjugation, we have
\[
g^*A_i
=\ii\bigl[(\lambda\alpha_i)\wedge
(\overline{\lambda}\,\overline{\alpha_i})\bigr]
=\lambda\overline{\lambda}\,A_i=A_i,
\]
where $\lambda\overline{\lambda}=1$ because $\lambda$ is a root of unity.
The eigenvalues on $\beta_i$ and $\overline{\beta_i}$ are $\lambda^{-1}$
and $\lambda$, so $g^*B_i=\lambda^{-1}\lambda B_i=B_i$.
Thus $A_i$ and $B_i$ are real $(1,1)$-classes fixed by every $g\in G$
and hence belong to $V=H^{1,1}(X,\R)^G$. Exterior multiplication gives
\begin{equation}\label{eq:nullclasses}
Q(A_i,A_i)=Q(B_i,B_i)=0,\qquad
s_i:=Q(A_i,B_i)=\int_X U_i\wedge\overline{U_i}>0,
\end{equation}
where positivity follows from the hypothesis $U_i=\alpha_i\wedge\beta_i\ne0$:
the integrand is nonnegative everywhere and positive on the nonempty open
set where $U_i$ does not vanish.

For $i\ne j$, we claim that the $2\times2$ matrix of intersections
between $(A_i,B_i)$ and $(A_j,B_j)$ has determinant $0$.
To see this, put $z_{ij}=\int_X U_i\wedge\overline{U_j}$.
Consider the two pairs
\[
\bigl([\alpha_i\wedge\overline{\alpha_j}],
      [\beta_j\wedge\overline{\beta_i}]\bigr),\qquad
\bigl([\alpha_i\wedge\overline{\beta_j}],
      [\alpha_j\wedge\overline{\beta_i}]\bigr).
\]
The two classes in the first pair belong to
$H^{1,1}(X)_{\chi_i-\chi_j}$, and the two classes in the second
belong to $H^{1,1}(X)_{\chi_i+\chi_j}$.
 Their Gram matrices for the Hermitian pairing
$(u,v)\mapsto\int_X u\wedge\overline v$ are
\begin{equation}\label{eq:grammatrices}
\begin{pmatrix}
-Q(A_i,A_j)&z_{ij}\\
\overline{z_{ij}}&-Q(B_i,B_j)
\end{pmatrix},\qquad
\begin{pmatrix}
-Q(A_i,B_j)&-z_{ij}\\
-\overline{z_{ij}}&-Q(B_i,A_j)
\end{pmatrix}.
\end{equation}
Both have rank at most 1 by \eqref{eq:mixedbound}. Their vanishing
determinants therefore give the single identity
\begin{equation}\label{eq:mixedidentity}
Q(A_i,A_j)Q(B_i,B_j)=|z_{ij}|^2=Q(A_i,B_j)Q(B_i,A_j).
\end{equation}

To bound the number of pairs, regard each $(A_i,B_i)$ as a single
vector $w_i=A_i\wedge B_i$ in the exterior square $\bigwedge^2_\R V$.
This exterior product is taken in the real vector space $V$.
The cup product $A_i\cup B_i$ on $X$ gives only the intersection number
$s_i$ after integration. The determinant of the $Q$-pairings is alternating
in each pair of vectors, so it defines a bilinear form on the exterior square:
\begin{equation}\label{eq:exteriorform}
Q_2(A\wedge B,C\wedge D)
=\det\begin{pmatrix}
Q(A,C)&Q(A,D)\\
Q(B,C)&Q(B,D)
\end{pmatrix}.
\end{equation}
This form is symmetric because interchanging the pairs transposes the matrix.
Equation~\eqref{eq:mixedidentity} now says that $Q_2(w_i,w_j)=0$ for $i\ne j$,
while \eqref{eq:nullclasses} gives
\[
Q_2(w_i,w_i)=\det\begin{pmatrix}0&s_i\\s_i&0\end{pmatrix}=-s_i^2<0.
\]
Thus the Gram matrix of $w_1,\ldots,w_m$ is diagonal with entries $-s_i^2$,
so these vectors span an $m$-dimensional negative definite subspace.
To bound $m$, we compute the negative index of $Q_2$ from the signature
$(1,t)$ of $Q$.
Choose a $Q$-orthogonal basis $e_0,\ldots,e_t$ with
$Q(e_0,e_0)=1$ and $Q(e_j,e_j)=-1$ for $j>0$.
The induced basis of $\bigwedge^2_\R V$ is $Q_2$-orthogonal:
the $t$ vectors $e_0\wedge e_j$ have square $-1$, while the vectors
$e_i\wedge e_j$ with $1\le i<j\le t$ have square $1$.
The negative index of $Q_2$ is therefore $t$, so $m\le t=\dim_\R V-1$.
\end{proof}

\subsection{The non-existence results}\label{sec:finish}

\begin{theorem}\label{thm:main}
There is no complex quintic surface in $\PP^3$ whose singular locus
consists of exactly 20 ordinary cusps.
\end{theorem}

\begin{proof}
Assume that such a 20-cuspidal quintic surface $Y$ exists.
By Proposition~\ref{prop:kfive}, its cusp code has
dimension 5. Proposition~\ref{prop:enumerator} supplies 10
representatives $r_1,\ldots,r_{10}$ of weight-18 words up to sign,
with all mixed sums and differences of weight 12.

Construct the smooth full cover $X$ of Subsection~\ref{sec:covers}.
By Proposition~\ref{prop:nonzero}, for each $i$ there are forms
$\alpha_i,\beta_i$ of characters $r_i,-r_i$ with
$\alpha_i\wedge\beta_i\ne0$. Proposition~\ref{prop:charspaces} gives
\[
\dim H^{1,1}(X)_{r_i-r_j}=
\dim H^{1,1}(X)_{r_i+r_j}=1\quad(i\ne j),
\]
and the invariant real $(1,1)$ space has signature $(1,4)$.
Theorem~\ref{thm:pair} now yields $10\le4$, a contradiction.
\end{proof}

\begin{corollary}\label{cor:canonical-nine}
There is no smooth complex surface of general type whose canonical map
is of degree 9 onto a surface of positive geometric genus.
\end{corollary}
\begin{proof}
We may pass to the minimal model without changing the canonical image
or the degree of the map. Theorem~\ref{thm:degree-nine} would make the
image a normal quintic
with exactly 20 ordinary cusps and no other singularities,
contradicting Theorem~\ref{thm:main}.
\end{proof}

\bibliographystyle{amsplain}
\bibliography{twenty_cusps}

\vspace{1cm}
\noindent\begin{minipage}{\textwidth}
\noindent Lev Borisov \vspace{0.1cm}
\\Department of Mathematics, Rutgers University
\\Piscataway, NJ 08854
\\ \verb|borisov@math.rutgers.edu|

\vspace{1cm}

\noindent Carlos Rito
\vspace{0.1cm}
\\ Centro de Matem\'atica, Universidade do Minho -- Polo CMAT-UTAD
\vspace{0.1cm}
\\ Universidade de Tr\'as-os-Montes e Alto Douro, UTAD
\\ Quinta de Prados
\\ 5000-801 Vila Real, Portugal
\vspace{0.1cm}
\\ \url{https://www.utad.pt}, \texttt{crito@utad.pt}

\vspace{1cm}

\noindent Xavier Roulleau \vspace{0.1cm}
\\Universit\'e d'Angers,
\\CNRS, LAREMA, SFR MATHSTIC,
\\F-49000 Angers, France
\\ \verb|xavier.roulleau@univ-angers.fr|
\end{minipage}
\end{document}